\documentclass[10pt,english,pre,nofootinbib,reprint,showpacs,showkeys,preprintnumbers,floatfix,longbibliography]{revtex4-1}
\usepackage[latin9]{inputenc}
\usepackage[a4paper]{geometry}

\usepackage{xcolor}
\usepackage{pdfcolmk}
\usepackage{graphicx}
\usepackage{babel}
\usepackage{mathrsfs}
\usepackage{amsmath}
\usepackage{amssymb}
\usepackage{esint}
\usepackage{subcaption}
\PassOptionsToPackage{normalem}{ulem}
\usepackage{ulem}
\usepackage{dsfont} 
\usepackage[unicode=true,pdfusetitle,
 bookmarks=true,bookmarksnumbered=false,bookmarksopen=false,
 breaklinks=false,pdfborder={0 0 1},backref=false,colorlinks=false]
 {hyperref}

\makeatletter

 \@ifundefined{textcolor}{}
 {%
   \definecolor{BLACK}{gray}{0}
   \definecolor{WHITE}{gray}{1}
   \definecolor{RED}{rgb}{1,0,0}
   \definecolor{GREEN}{rgb}{0,1,0}
   \definecolor{BLUE}{rgb}{0,0,1}
   \definecolor{CYAN}{cmyk}{1,0,0,0}
   \definecolor{MAGENTA}{cmyk}{0,1,0,0}
   \definecolor{YELLOW}{cmyk}{0,0,1,0}
 }

\usepackage{aas_macros}
\usepackage{bbold}

\let\textquotedbl="

\makeatother

\newtheorem{criterion}{Criterion}
\begin{document}

\title{Optimal Belief Approximation}

\author{Reimar H. Leike, Torsten A. Enßlin}

\affiliation{{\small{}Max-Planck-Institut für Astrophysik, Karl-Schwarzschildstr.~1,
85748 Garching, Germany}\\
Ludwig-Maximilians-Universität München, Geschwister-Scholl-Platz{\small{}~}1,
80539 Munich, Germany}
\begin{abstract}
In Bayesian statistics probability distributions express beliefs. However, for many problems
the beliefs cannot be computed analytically and approximations of
beliefs are needed. We seek a loss function that quantifies how ``embarrassing'' it is to communicate a given approximation. 
We reproduce and discuss an old proof showing that there is only one ranking under the requirements that (1) the best ranked approximation is the 
non-approximated belief and
(2) that the ranking judges approximations only by their predictions for actual outcomes.
The loss function that is obtained in the derivation is equal to the Kullback-Leibler divergence when normalized.
This loss function is frequently used in the literature. 
However, there seems to
be confusion about the correct order in which its functional arguments---the~approximated and
non-approximated beliefs---should be used. 
The correct order ensures that the recipient of a communication is only deprived of the minimal amount of information.
We hope that the elementary derivation settles the apparent confusion. 
For example when approximating beliefs with Gaussian
distributions the optimal approximation is given by moment matching.
This is in contrast to many suggested computational schemes.
\end{abstract}

\maketitle

\section{Introduction}

In Bayesian statistics, probabilities are interpreted as degrees of belief. For any set of mutually exclusive and exhaustive events, one expresses the state 
of knowledge as a probability distribution over that set. The probability of an event then describes the personal confidence that this event will happen or has
happened. As a consequence, probabilities are subjective properties reflecting the amount of knowledge an observer has about the events; a different observer
might know which event happened and assign different probabilities. If an observer gains information, she updates the probabilities she had assigned before.

If the set of possible mutually exclusive and exhaustive events is infinite, it is generally impossible to store all entries of the corresponding probability distribution on a
computer or communicate it through a channel with finite bandwidth.
One therefore needs to approximate the probability distribution which describes one's belief. Given a limited set $X$ of approximative
beliefs $q(s)$ on a quantity $s$, what is the best belief to approximate the actual belief as expressed by the probability $p(s)$?

In the literature, it is sometimes claimed that the best approximation is given by the $q\in X$ that minimizes the Kullback--Leibler divergence (``approximation'' KL) \citep{Kullback1951}
\begin{align}
\text{KL}(p,q)\label{eq:KL-approximation}= \sum_sp(s)\,\text{ln}\left(\frac{p(s)}{q(s)}\right)
\end{align}
where $q$ is the approximation and $p$ is the real belief. We refer to this functional as ``approximation'' KL 
to emphasize its role in approximation, which will be derived in the course of this paper and to distinguish it from the same functional, with $q$ being a prior belief
and $p$ being the posterior belief to which this KL is minimized in inference. 
We will refer to the functional with $q$ being the input and $p$ obtained through minimization
~as ``inference KL''.
In Equation (\ref{eq:KL-approximation}), minimization is done with respect to its second argument. The derivation of this particular functional form varies from field to~field.

For example, in coding theory, one tries to minimize the amount of bandwidth needed to transmit a message. Given a prior $q$ over the symbols that the message
consists of, an optimal scheme can be derived. The approximation KL gives the expected amount of extra bits needed to transmit such a message if the symbols are actually drawn from the probability distribution
$p$ instead of $q$ \citep{cover2006elements}. 
If we know that $p$ is the real probability distribution, 
the best approximative
probability distribution $q\in X$ to base a coding on is therefore the one minimizing the approximation KL.
However, it is  not clear that minimizing the amount of bits transferred is the best or even only
measure expressing how good such an approximation is in general.

In machine learning and deep learning, neural networks are trained to understand abstract data $d$; for example, to assign a label $s$ to it. This task can be viewed as fitting an approximative probability distribution $q(s|d)$ to a true generating
probability distribution $p(s|d)$. For this, the~approximative probability distribution is parametrized (to a neural network) and then matched to the true, generating probability distribution
using a loss function and samples. The most frequently used loss function is the cross entropy, which is equivalent to the approximation KL. The reason to use
this form is often either inspired from coding theory, or by experimental experience \citep{Goodfellow-et-al-2016}.

Another argument for minimizing the approximation KL is given in Chapter 13 of Reference \citep{opper2001advanced}, where it is claimed that
this yields the maximum likelihood estimation to $p(s)$ among the probability distributions in $X$ and that it gives an unbiased and unique approximation.
Interchanging the arguments of the Kullback--Leibler divergence (the inference KL used in variational Bayes) generally leads to a biased estimate and does not necessarily yield a unique result. 
These arguments undoubtedly give evidence for why minimizing the approximation KL gives a good estimate. However, this does not exclude all other methods. 
Having an unbiased estimate refers to getting the right mean. In our picture, this is a result of optimal approximation and not a requirement for optimality.
Additionally, this~result was derived with the help of information geometry, whose applicability to non-local problems is criticized, for example, in
References\,\citep{skilling2014critique, skilling2015failures}.

Contrary to the evidence for minimizing the approximation KL, we find many examples where an approximation is made by minimizing other functionals;
for example, minimizing the inference KL (e.g.,\,\citep{fox2012tutorial, 2010PhRvE..82e1112E, 2016arXiv160100670B, 2014arXiv1408.1920P, 2013arXiv1310.7845P, 2013PhRvE..87a3308E, 2008PhyA..387.6759T}).
For many but not all of them, this is because minimizing the approximation KL is not feasible in practice in their case due to the real distribution $p$ not being~accessible.

In this paper, we seek to bring together the different motivations and give a full and consistent picture.
The proof we present here is not new; it goes back to \citep{bernardo1979expected}, where there is an exact mathematical derivation for probability densities 
analogously to our derivation.
However, there are earlier publications dealing with the discrete case \citep{aczel1967remarks, mccarthy1956measures, good1952rational}.
Although this proof dates back at least 40 years, its implication on approximating beliefs seem to be quite unknown---especially in the community of physicists applying Bayesian methods. 
In this paper, we reproduce a slightly modified version of this proof, give the result a new interpretation and
add further justification for the prerequisites used, laying emphasis on why one has to accept the axioms necessary for the derivation if one is a Bayesian.
We 
provide argumentation for why believe approximation is an important and omnipresent topic. 

We lay the emphasis of this paper more on interpretation of results and justification of prerequisites, and thus present an easy version of the proof
where the loss function is assumed to be differentiable.
The proof can however be extended to the general case of non-differentiable loss \citep{harremoes2017divergence}.
The argument we reproduce gives evidence that minimizing the approximation KL is the best approximation in theory. This argument
does not rest on information geometry nor is it restricted to coding theory.
By imposing two consistency requirements, one is able to exclude all ranking functions with the exception of
~one for
ranking the approximative probability distributions $q\in X$.
For this, one employs the principle of loss functions \citep{cramer1930mathematical}, also called cost functions, regret~functions, 
(or with flipped sign, utility functions or score functions) and shows that the unique loss function for ranking 
approximated probability distributions is the approximation KL. 
For us, a ranking is a total order indicating preference, whereas a loss is 
a map to $\mathbb{R}$, which induces a ranking but additionally gives an absolute scale to compare preferences.
The presented axiomatic derivation does not give rise to any new method, but 
it enables a simple checking for whether a certain
approximation is most optimally done through the approximation KL. 
There are many other examples of axiomatic derivations seeking to support information theory on a fundamental level. Some notable examples are
Cox derivation \citep{Cox1946} of Bayesian probability theory as a unique extension of Boolean algebra as well as a scientific discussion on the 
maximum entropy principle \citep{2003prth.book.....J, skilling1988axioms, 2004AIPC..707...75C}, establishing
the inference KL as unique inference tool (which gave rise to the naming convention in this paper). Most~of these arguments rely on page-long proofs
to arrive at the Kullback--Leibler divergence. The proof that is sketched in this paper is only a few lines long, but nonetheless
standard literature for axiomatic derivation in Bayesianism does not cite this ``easy'' derivation (e.g., the influential Reference\,\citep{2003prth.book.....J}). 
As~already discussed, approximation is an
important and unavoidable part of information theory, and with the axiomatic derivation presented here we seek to provide orientation to
scientists searching for a way to approximate probability~distributions.

In Section \ref{sec:Introduction-to-Loss}, we introduce the concept of loss
functions, which is used in Section \ref{sec:Deriving-the-Unique} to 
define an optimal scheme for approximating probability distributions that express beliefs. We briefly discuss the relevance of our derivations
for the scientific community in Section \ref{sec:Discussion}.
We conclude in Section \ref{sec:Conclusion}. 

\section{Loss Functions\label{sec:Introduction-to-Loss}}

The idea to evaluate predictions based on loss functions dates back 70 years, and
was first introduced by Brier \citep{brier1950verification}.
We explain loss functions by the means of parameter estimation. Imagine that one would like to give an estimate of a parameter $s$ that is not known,
which value of $s$ should be taken as an estimate? 
One way to answer this question is by using loss functions.
For this note that $p(s)$ is now formally a probability measure, however we choose to write
$
\int \text{d}s p(s)
$
instead of
$
\int \text{d}p(s)
$
as if $p(s)$ would be a probability density.
A loss function in the setting of parameter estimation
is a function that takes an estimate $\sigma$ for $s$ and quantifies how ``embarrassing'' this estimate is if $s=s_0$ turns out to be the case:
\begin{align*}
 \mathscr{L}(\sigma,s_0)
\end{align*}

The expected embarrassment can be computed by using the knowledge $p(s)$ about $s$:
\begin{align*}
\left<\mathscr{L}(\sigma,s_0)\right>_{p(s_0)}=\int \text{d}{s_0}\,\mathscr{L}(\sigma,s_0)\, p(s_0)
\end{align*}

The next step is to take the estimate $\sigma$ that minimizes the expected embarrassment; that is,
the~expectation value of the loss function.
For different loss functions, one arrives at different recipes for how to extract an estimate $\sigma$
from the belief $p(s)$; for example, for $s\in\mathbb{R}$:

\begin{align}
&\mathscr{L}(\sigma,s_0) = \nonumber\\
&=
\begin{cases}
-\delta(\sigma-s_0)\label{eq:delta-loss}  & \Rightarrow 
 \begin{array}{l}
\text{Take $\sigma$ such that}\\ 
 p(s)|_{s=\sigma} \text{ is maximal} 
 \end{array} \\  
\left|\sigma-s_0\right|  & \Rightarrow \text{ Take $\sigma$ to be the median} \\
(\sigma-s_0)^2  & \Rightarrow \text{ Take $\sigma$ to be the mean} \\
\end{cases}
\end{align}

In the context of parameter estimation, there is no general loss function that one should take. In~many scientific applications, the third option is favored, 
but different situations might enforce different loss functions.
In the context of probability distributions, one has a mathematical structure available to guide the choice. In this context, one can 
restrict the possibilities by requiring consistent loss functions. 

\section{The Unique Loss Function\label{sec:Deriving-the-Unique}}

How embarrassing is it to approximate a probability distribution by $q(s)$ even though it is actually $p(s)$? We quantify the embarrassment in a loss function
\begin{align}
\mathscr{L}\left(\frac{q}{m},s_0\right)  \label{eq:loss-def}
\end{align}
which says how embarrassing it is to tell someone $q(s)$ is one's belief about $s$ in the event that later $s$ is measured to be $s_0$. Here $m$ is introduced as reference measure to make 
$\mathscr{L}$ coordinate independent. For a finite set coordinate independence is trivially fulfilled and it might seem that
having a reference measure $m$ is superficial. Note however, that it is a sensible additional requirement to have the quantification be invariant under splitting of events, 
i.e. mapping to a bigger set where two now distinguishable
events represent one former large event. The quotient $\frac{q}{m}$ is invariant under such splitting of events, whereas $q$ itself is not. 
The reference measure $m$ can be any measure such that $q$ is absolutely continuous with respect to $m$.

Note further that we restrict ourselves to the 
case that we get to know the exact value of $s$. This does not make our approach less general;
imagining that we would instead take a more general loss $L(q,\tilde{q}(s))$ where $\tilde{q}$ is the knowledge about $s$ at some
later point, then we may define $\mathscr{L}\left(\frac{q}{m},s_0\right) = L(q,\delta_{ss_0})$ with $\delta$ denoting the Kronecker or
Dirac delta function, and thus restrict ourselves again to the case of exact knowledge.
This line of reasoning was spelled out in detail by John Skilling \citep{skilling1989classic}:
\begin{quote}
``If there are general theories, then they must apply to special cases''.
\end{quote}

To decide which belief to tell someone, we look at the expected loss
\begin{align}
\left<\mathscr{L}\left(\frac{q}{m},s_0\right) \right>_{p(s_0)}=\int\text{d}s_0\,\mathscr{L}\left(\frac{q}{m},s_0\right) \, p(s_0)
\end{align}
and try to find a $q\in X$ that minimizes this expected loss.
To sum up, if we are given a loss function, we have a recipe for
how to  optimally approximate the belief.
Which loss functions are sensible, though? We enforce two criteria that a good loss function should satisfy. 
\begin{criterion}\emph{(Locality)}
If $s=s_0$ turned out to be the case, $\mathscr{L}$ only depends on the prediction $q$ actually makes about $s_0$:
\begin{align}
\mathscr{L}\left(\frac{q}{m},s_0\right)=\mathscr{L}\left(\frac{q(s_0)}{m(s_0)}\right)\label{eq:loss-local}
\end{align}
\end{criterion}
Note that we make an abuse of notation here, denoting the function on both sides of the equation by the same symbol.
The criterion is called locality because it demands that the functional of $\frac{q}{m}$ should be
evaluated locally for every $s_0$. It also forbids a direct dependence of the loss $\mathscr{L}$ on $s_0$ which
excludes losses that are a priori biased towards certain outcomes $s_0$.

This form of locality is an intrinsically Bayesian property.
Consider a situation where one wants to decide which of two rivaling hypotheses to believe. In order to distinguish them, some data $d$ are measured.
To update the prior using Bayes theorem, one only needs to know how probable the \emph{measured} data $d$ are given each hypothesis, not how probable
other \emph{possible} data $\tilde{d}\neq d$ that were not measured are. This might seem intuitive, but there exist hypothesis decision methods (not necessarily based on loss functions)
that do not fulfill this property. For example, the non-Bayesian $p$-value depends mostly on data that were not measured (all the data that are at least as ``extreme'' as the measured data).
Thus, it is a property of Bayesian reasoning to judge predictions only by using what was predicted about things that were measured.


The second criterion is even more natural. If one is 
not restricted in what can be told to others, then the best thing should be to tell them the actual belief $p$.
\begin{criterion}\emph{(Optimality of the actual belief, properness.)}
Let $X$ be the set of all probability distributions over $s$. For all $p$ and all $m$, the probability distribution $q\in X$ with minimal expected loss is
$q=p$:
\begin{align}
0 = \left(\frac{\partial}{\partial q(s)}\left<\mathscr{L}\left(\frac{q}{m},s_0\right)\right>_p\label{eq:posterior-optimal}\right)_{q=p}
\end{align}
\end{criterion}

The last criterion is also referred to as a proper loss (score) function in the literature; see  Reference \citep{gneiting2007strictly} for a mathematical overview of different
proper scoring rules. Our version of this property is slightly modified to the version that is found in the literature as we
demand this optimum to be obtained independently of a reference measure $m$. There is a fundamental Bayesian desiderata stating
that ``If there are multiple ways to arrive at a solutions, then they must agree.'' We'd like to justify why this is a property that is
absolutely important. If one uses statistics as a tool to answer some question, then if that answer 
would depend on how statistics is applied, then this statistic itself is inconsistent.
In our case, where the defined loss function is dependent on an arbitrary reference measure $m$, the result is thus forced to be independent
of that $m$. 

Note furthermore that although intuitively we want the global optimum to be at the actual belief $p$ (referred to as strictly proper in the literature), 
mathematically we only need it to be an extreme value for our derivation.

Having these two consistency requirements fixed, we derive which kind of consistent loss functions are possible. We~insert Equation~(\ref{eq:loss-local}) into Equation  (\ref{eq:posterior-optimal}), expand the domain of the loss function to not necessarily normalized positive vectors $q(s)$ but 
introduce $\lambda$ as a Lagrange multiplier to account for the fact that we minimize under the constraint of
normalization. We compute
\begin{equation}
\begin{array}{lllll}
0 &=& \left(\frac{\partial}{\partial q(s)}\int \text{d}{s_0}\,\left(\mathscr{L}\left(\frac{q(s_0)}{m(s_0)}\right)\, p(s_0)+\lambda\, q(s_0)\right)\right)_{q=p}\\
&=& \int \text{d}{s_0}\,\partial\mathscr{L}\left(\frac{p(s_0)}{m(s_0)}\right)\frac{\delta(s-s_0)}{m(s_0)}\,p(s_0)+\lambda\,\delta(s-s_0)\\
&=& \partial\mathscr{L}\left(\frac{p(s)}{m(s)}\right)\frac{p(s)}{m(s)}+\lambda\\
&\Rightarrow& \partial\mathscr{L}\left(\frac{p(s)}{m(s)}\right) = -\frac{\lambda m(s)}{p(s)}\label{eq:loss-with-quotient}
\end{array}
\end{equation}
Here $\partial\mathscr{L}$ denotes the derivative of $\mathscr{L}$.
In the next step we substitute $x:=\frac{p(s)}{m(s)}$ for the quotient. Note that Equation~(\ref{eq:loss-with-quotient}) holds for all positive real values of $x\in \mathbb{R}^+$ since
the requested measure independence of the resulting ranking permits to insert any measure $m$. We then obtain
\begin{equation}
\begin{array}{lll}
\partial\mathscr{L}\left(x\right) &=& -\frac{\lambda}{x}\\
\Rightarrow \mathscr{L}\left(x\right) &=& -C\,\text{ln}\left(x\right)+D
\end{array}
\end{equation}
where $C>0$ and $D$ are constants with respect to $q$.
Note that through the two consistency requirements one is able to completely fix what~the loss function is. In the original proof of \citep{bernardo1979expected}, 
there arises an additional possibility for $\mathscr{L}$ 
if the sample space consists of 2 elements. In that case, the locality axiom, as it is used in the literature, does not constrain $\mathscr{L}$ at all. 
In our case, where we introduced $m$ as a reference measure,
we are able to exclude that possibility.
Note that the constants $C$ and $D$ are irrelevant for determining the optimal approximation as they do not affect where the minimum of the loss function is.

To sum up our result if one is restricted to the closed set $X$ of probability distributions, one should take that $q\in X$ that minimizes
\begin{align}
\left<\mathscr{L}\left(\frac{q}{m},s_0\right)\right>_{p(s_0)}=-\int\text{d}s_0\,p(s_0)\,\text{ln}\left(\frac{q(s_0)}{m(s_0)}\right)\label{cross-entropy}
\end{align}
in order to obtain the optimal approximation, where it is not important what $m$ is used.

If one takes $m=1$, this loss is the cross entropy
$$
\left<-\text{ln}\left(q(s_0)\right)\right>_{p(s_0)}\ .
$$
If one desires a rating of how good an approximation is, and not only a ranking which approximation is best, one could go one step further and enforce a third criterion:
\begin{criterion}\emph{(zero loss of the actual belief)}
For all $p$, the expected loss of the probability distribution $p$ is 0:
\begin{align}
0 = \left<\mathscr{L}\left(\frac{p}{m},s_0\right)\right>_p\label{eq:posterior-0}
\end{align}
\end{criterion}
This criterion trivially forces $m=p$ and makes the quantification unique while inducing the same ranking. Thus we arrive at the Kullback-Leibler divergence
$$
\text{KL}(p,q) = \int \text{d}s_0\, p(s_0)\text{ln}\left(\frac{p(s_0)}{q(s_0)}\right)
$$
as the optimal rating and ranking function.

To phrase the result in words, the optimal way to approximate 
the belief $p$
is such that given the approximated belief $q$, the amount of information $\text{KL}(p,q)$ that has to be obtained for someone who believes $q$ to arrive back at the actual belief $p$ is 
minimal. We should make it as easy as possible for someone who got an approximation $q$ to get to the correct belief $p$.
This sounds like a trivial statement, explaining why the approximation KL is already widely used for exactly that task.

\section{Discussion\label{sec:Discussion}}

We briefly discuss the implications of these results.

In comparison to Reference\,\citep{opper2001advanced, cover2006elements}, we presented another more elementary line of argumentation for the
claim that the approximation KL is the correct ranking function for approximating which holds in a more general setting. 

Other works that base their results on minimizing the inference KL ($\text{KL}(q,p)$) for belief approximation are not optimal with respect to the ranking function we derived in
Section\,\ref{sec:Deriving-the-Unique}. 
One~reason for preferring the for this purpose non-optimal inference KL is that it is computationally feasible for many applications, in contrast to the optimal 
approximation. As long as the optimal scheme is not computationally accessible, this argument has its merits.

Another reason for minimizing the inference KL for approximation that is often cited (\mbox{e.g., \citep{bishop2007pattern}})
is that it gives a lower bound to the log-likelihood
\begin{align}
\text{ln}(p(d|s)) = \left<\text{ln}\left(\frac{p(d,s)}{q(s)}\right)\right>_{q(s)} + \text{KL}(q,p)
\end{align}
which for example gives rise to the expectation maximization (EM-
) algorithm \citep{dempster1977maximum}. However, the method only gives rise to maximum a posteriori or maximum likelihood
solutions, which corresponds to optimizing the $\delta$-loss of Equation\,(\ref{eq:delta-loss}).

\begin{figure}[h]
     \centering
   \includegraphics[width=.5\textwidth]{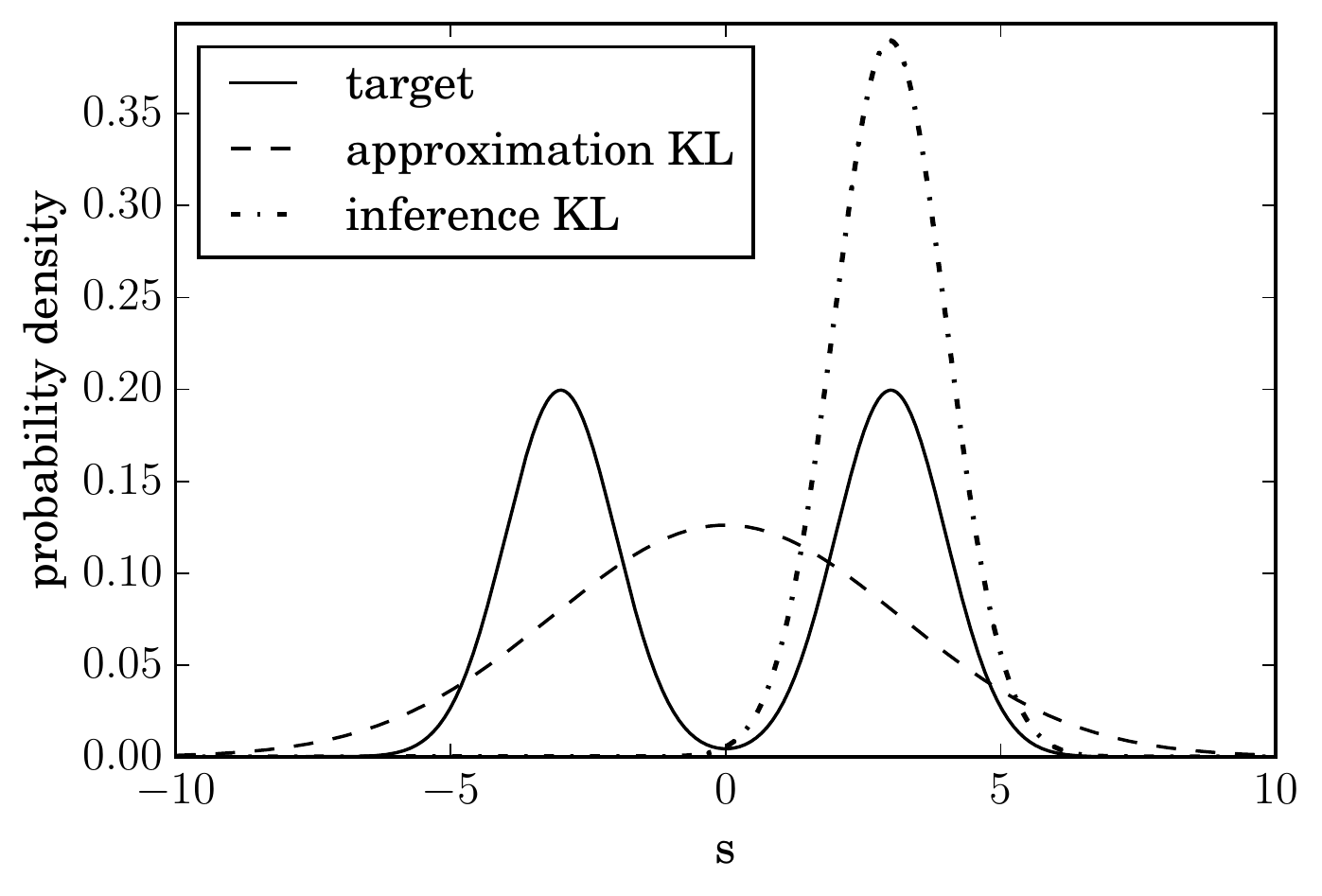}
   \vspace{-6pt}

   \caption{
   Results of approximating a target distribution in $s$ with a Gaussian distribution. KL: Kullback--Leibler.}
   \label{fig:example-KL}
\end{figure}

In Reference\,\citep{2013arXiv1310.7845P}, it is claimed that minimizing the inference KL yields more desirable results since
for multi-modal distributions, individual modes can be fitted with a mono-modal distribution such as a Gaussian distribution, whereas the resulting distribution has a very large variance when minimizing the approximation KL
to account for all modes. In Figure\,\ref{fig:example-KL} there is an example of this behavior. 
The~true distribution of the quantity $s$ is taken to be a mixture of two standard Gaussians with means $\pm3$. 
It is approximated with one Gaussian distribution by using the approximation KL and the inference KL. When using the approximation KL, the resulting distribution has a large variance to cover both peaks. 
Minimizing the inference KL leads to a sharply peaked approximation around one peak. 
A user of this method might be very confident that the value of $s$ must be near $3$, even
though the result is heavily dependent on the initial condition of the minimization and could have become peaked around $-3$ just as well.

We find that fitting a multi-modal distribution with a mono-modal one will yield suboptimal results irrespective of the fitting scheme. An approximation should always
have the goal to be close to the target that is being approximated. If it is already apparent that this goal cannot be achieved, it is recommended to rethink the set of approximative
distributions and not dwell on the algorithm used for approximation.

In Reference\,\citep{2013PhRvE..87a3308E} an approximative simulation scheme called information field dynamics is described.
There, a Gaussian distribution $q$ is matched to a time-evolved version $U(p)$ of a Gaussian distribution $p$. 
This matching is done by minimizing the inference KL. 
In this particular case (at least for information preserving dynamics), the matching can be made optimal without making the algorithm more complicated. 
Since for information preserving dynamics time evolution is
just a change of coordinates and the Kullback--Leibler divergence is invariant under such transformations, one can instead match the Gaussian distribution $p$ and $U^{-1}(q)$ by 
minimizing $\text{KL}(p,U^{-1}(q)) = \text{KL}(U(p),q)$, which is just as difficult in terms of computation.
 \vspace{-6pt}

In Reference\,\citep{2008PhyA..387.6759T} it is claimed that the inference KL yields an optimal approximation scheme fulfilling certain axioms. 
This result is the exact opposite of our result. This disagreement is due to an assumed consistency of approximations. In Reference\,\citep{2008PhyA..387.6759T}, 
further approximations are forced to be consistent with earlier approximations; i.e., if one does two approximations, one gets the same result as with one joint approximation.
Due to this requirement, the derived functional cannot satisfy some of the axioms that we used. In our picture, it is better to do one large approximation
instead of many small approximations. This is in accordance to the behavior of other approximations. For example, when step-wise rounding the real number $1.49$, one gets
$2$ if it is first rounded to one decimal and then to integer precision compared to being rounded to integer precision directly where one gets $1$. If~information is lost due to approximation, it is natural for further approximations to be less precise than if one were to approximate in one go.

There also exist cases where we could not find any comments explaining why the arguments of the Kullback--Leibler divergence are in that particular order. 
In general, it would be desirable that authors provide a short argumentation for why 
they choose a particular order of the arguments of the KL divergence.

\section{Conclusions\label{sec:Conclusion}}

Using the two elementary consistency requirements on locality and optimality, as expressed by~Equations\ (\ref{eq:loss-local}) and (\ref{eq:posterior-optimal}), respectively, we have shown that there is only 
one ranking function that ranks how good an approximation of a belief is, analogously to Reference\,\citep{bernardo1979expected}.
By minimizing $\text{KL}(p,q)$ with respect to its second argument $q\in X$, one gets the best approximation to $p$. This is claimed at several points in the literature.
Nevertheless, we found sources where other functionals were minimized in order to obtain an approximation.
This confusion is probably due to the fact that for the slightly different task of updating a belief $q$ under new constraints, $\text{KL}(p,q)$ has to be minimized with respect
to $p$, its first argument \citep{csiszar1991least, 2014arXiv1412.5644C}. We do not claim that any of the direction of Kullback--Leibler divergence are wrong by themselves, but one should
be careful of when to use which.

We hope that for the case of approximating a probability distribution $p$ by another $q$ we have given convincing and  conclusive arguments for why this should be done by minimizing
$\text{KL}(p,q)$ with respect to $q$, if this is feasible.

\section{Acknowledgements}
We would like to thank A. Caticha, J. Skilling, V. B\"ohm, J. Knollm\"uller, N. Porqueres, and M. Greiner and six anonymous referees for the discussions and their valuable comments on the manuscript.

\bibliographystyle{apsrev4-1}
\bibliography{ift}
\end{document}